\documentclass{amsart}
\usepackage{amsfonts}
\usepackage{amsmath}
\usepackage{amssymb}
\usepackage{graphicx}
\newtheorem{theorem}{Theorem}

\begin{document}

\markboth{A.M. Baxter}{Journal of Difference Equations and Applications}

\title{Applying the Cluster Method to Count Occurrences of Generalized Permutation Patterns}
\author{Andrew M. Baxter}

\begin{abstract}
We apply ideas from the cluster method to $q$-count the permutations of a multiset according to the number of occurrences of certain generalized patterns, as defined by Babson and Steingr\'{i}msson.  In particular, we consider those patterns with three letters and one internal dash, as well as permutation statistics composed of counting the number of occurrences of multisets of such patterns.  Counting is done via recurrences which simplify in the case of permutations.  A collection of Maple procedures implementing these recurrences accompanies the article.
\end{abstract}


\maketitle

\section{Introduction}

In its most basic form the cluster method, as described in \cite{NoonZeil}, counts the number of words of length $\ell$ in a given alphabet which avoid a certain set of forbidden subwords (i.e. contiguous blocks of letters).  In that same article, Noonan and Zeilberger extend the method to count the number of words of length $\ell$ which contain a given number of those forbidden subwords.  The permutation statistics $des$, $inv$, and $maj$ were implicitly treated with this method in \cite{Zeil80}.  Here we consider distributions of generalized permutation patterns of length 3, as introduced by Babson and Steingr\'{i}msson \cite{BabStein}.  Burstein and Mansour \cite{BurstMans} lists the number of words of length $\ell$ in the alphabet $\{1,..,k\}$ avoiding a given generalized pattern for each pattern of length 3.  Claesson and Mansour \cite{ClaesMans} provide a recurrence to compute the number of permutations with $r$ occurences of a generalized pattern of length 3 with one internal dash.  In this paper, we derive recurrences to compute (via a generating function) the number of permutations of a multiset with $r$ occurences of a generalized pattern of length 3 with one internal dash.  The special case of permutations of $\{1,\ldots,n\}$ allows for some more specialized recurrences.  The method can simultaneously compute the multivariate distribution of different patterns, yielding recurrences to compute certain pattern-based permutation (multi)statistics.

In section 2 we outline the conventions regarding operations on words and recall definitions of generalized patterns.  Section 3 describes how we apply the cluster method to derive recurrences for the distributions of generalized patterns on multiset permutations.  Section 4 considers the special case of distributions of patterns over permutations in $S_n$, adapting the recurrences found in Section 3.  Section 5 details the accompanying Maple program which provides several tools for the reader to verify and to explore the conclusions discussed in this paper.  A summary of results and suggestions of future directions can be found in the final section.

\section{Preliminaries}

\subsection{Operations on Words and Notational Conventions}

We will consider words $w=w_1 w_2 \cdots w_\ell$ in the alphabet $A=[n]=\{1,2,\ldots,n\}$.  Borrowing notation from the theory of formal languages, let $A^*$ denote the set of all (finite) words in the alphabet $A$.  We will use $\ell(w)$ or simply $\ell$ to represent the length of word $w$ and $n$ to denote the alphabet size.  In the case of permutations we use $n$ to represent both.  For words $u=u_1 \cdots u_k$ and $v=v_1 \cdots v_\ell$, we denote the concatenation as $uv = u_1 \cdots u_k v_1 \cdots v_\ell$.  We use $\nu_i(w)$ to denote the multiplicity of $i$ in $w$, and in generating functions it will be useful to consider $x^w := \prod_{i\geq 1} x_i^{\nu_i(w)}$, where the $x_i$ are commuting indeterminates.  Note that $x^{uv}=x^u x^v$. 

We may also define the reversal and complement operations on words in $[n]$.  Let $w^r = w_\ell w_{\ell-1} \cdots w_1$ and $w^c = (n+1-w_1) (n+1-w_2) \cdots (n+1-w_\ell)$.  If we are given an alphabet vector $\mathbf{m}=(m_1, \ldots, m_n)$ and consider the set of words such that $\nu_i (w) = m_i$ for each $i$, then reversal is a bijection.  Complementation is a bijection if and only if $m_1=m_n, m_2=m_{n-1}, \ldots, m_{\lfloor n/2 \rfloor}=m_{\lceil n/2 \rceil}$.

Since we count the descents of a word, it will be useful to define the set $Des(w) := \{i: w_i > w_{i+1}\}$.  This should not be confused with $des(w):=|Des(w)|$.

\subsection{Generalized Pattern Functions}

Babson and Steingr\'{i}msson introduced generalized permutation patterns in \cite{BabStein}, giving definitions which apply equally well to words (i.e. permutations with repeated letters).  For simplicity we work with words in the alphabet $[n]=\{1,2,\ldots, n\}$.  A pattern is written in the alphabet $\{a,b,c,\ldots\}$ where dashes may or may not separate successive letters.  The absence of a dash between two letters indicates that the corresponding letters in the word must be adjacent; otherwise corresponding letters may be arbitrarily far apart.  Further the corresponding letters of the word must be order-isomorphic to the pattern.  We will write the pattern as a function, where applying the function counts the number of occurrences of the pattern.  This is best understood through examples.  The pattern function $(a-bc)(w)$ for word $w=w_1 w_2 \cdots w_\ell$ is the number of subsequences $w_i w_j w_{j+1}$ such that $i<j$ and $w_i<w_j<w_{j+1}$.  Likewise $(abc)(w):=\left| \left\{ i: w_i<w_{i+1}<w_{i+2}  \right\}  \right|$ and $(ba)=| \{ i: w_{i}>w_{i+1} \} |$.  For words we have additional (nontrivial) pattern functions with repeated letters such as $(a-ab)(w)=| \{ (i,j): i<j,  w_i=w_j<w_{j+1} \} |$.  We will regularly use the symbols $\sigma$ and $\tau$ to represent pattern functions.  Treating these as functions allows multi-pattern functions like $\sigma+2\tau$, the number of occurences of $\sigma$ plus twice the number of occurences of $\tau$.  

In this paper, we will only explicitly consider those patterns with three letters and one internal dash between the first and second letters, that is those of the shape $(x-yz)$.  These are called \emph{type-(1,2)} by Mansour \cite{Mans}.  If $y>z$, then we call such patterns \emph{descent-based} and if $y<z$ they are called \emph{rise-based}.  Patterns of type-(2,1) are easily handled via the reversal operation, since $\sigma^r(w)=\sigma(w^r)$ where $(x-yz)^r = (zy-x)$.

We are interested in the distribution of patterns over sets of words.  For a pattern function $\sigma$ and set of words $W$, the \emph{distribution of $\sigma$ over $W$} is the generating function
\[
F_\sigma^W (q) = \sum_{w \in W} q^{\sigma(w)}.
\]
The $W$ superscript is usually omitted.  We will use ideas from the cluster method to derive recurrences which quickly compute $F_\sigma$ over the set of all words in the alphabet $[n]$, as well as the set of words with a prescribed number of each letter.  Two statistics $\sigma, \tau$ are \emph{equidistributed} over $W$ if $F_{\sigma}^W(q) = F_{\tau}^W(q)$.

\section{The Cluster Method for Recurrences}

We now adapt the cluster method to count the number of words with a given number of occurrences of a given pattern.  Instead of counting the number of occurences of a descent-based [resp. rise-based] pattern, we instead consider the problem of weight-counting the descents [resp. rises], where the descent-weight [resp. rise-weight] is determined by the pattern.  For example, the pattern $(a-cb)$ can be weight-counted,
\[
(a-cb)(w) = |\{(j,i): j<i, w_j < w_{i+1} < w_i\}| = \sum_{i \in Des(w)} | \{j: j<i, w_j < w_{i+1} \}|.
\]
Denote the descent-weight of $i\in Des(w)$ according to pattern $\sigma$ by $\sigma_d(w,i)$.  Hence $(a-cb)_{d}(w,i) = | \{j: j<i, w_j < w_{i+1} \}|$.  

Thus we may manipulate the sum in $F_\sigma$ as follows:

\begin{eqnarray*}
  F_{\sigma}(q,x)&=&\sum_{w\in A^*} x^w q^{\sigma(w)} \\
  &=& \sum_{w\in A^*} x^w \prod_{i \in Des(w)} q^{\sigma_d(i,w)} \\
  &=& \sum_{w\in A^*} x^w \prod_{i \in Des(w)} (1+(q^{\sigma_d(i,w)}-1)) \\
  &=& \sum_{w\in A^*} x^w \sum_{S\subseteq Des(w)} \prod_{i\in S} (q^{\sigma_d(i,w)}-1),
\end{eqnarray*}
  
\noindent where in the last step we expand the product.  In essence, this is the principle of Inclusion-Exclusion (alternately, M\"{o}bius inversion or sieve methods).  As in \cite{NoonZeil} we are counting pairs $(w,S)$ which can be considered \emph{marked words}, where some subset $S$ of the descents of $w$ are ``marked'' (as a teacher would mark a wrong answer).  Each marked descent $w_i w_{i+1}$ has weight $(q^{\sigma_d(w,i)}-1)$, and the weight of a marked word is the product of the weights of the marked descents.

We break here to provide a concrete example.  Consider the word $w=2637541$ and pattern $\sigma = (a-cb)$.  Then $\sigma_d(w,4)=2$, since $2-75$ and $3-75$ are both occurences of $(a-cb)$ involving the descent at position 4.  Similarly, $\sigma_d(w,2)=1$, $\sigma_d(w,5)=2$, and $\sigma_d(w,6)=0$.  Now consider the marked word $(w,\{2,5\})$, better viewed as $2\underline{63}7\underline{54}1$.  Only the descents 63 and 54 contribute, yielding the weight $(q^1 -1)(q^2 -1)$.  Similarly, the weight of $(w,\{4,5,6\})$ is $(q^2 -1)(q^2 -1)(q^0-1)$ and the weight of $(w,\emptyset)$ is 1 because of the empty product.  It is easily checked that the sum of the weights each of the $2^4$ marked versions of $w$ equals $q^{\sigma(w)}=q^5$

From the cluster method, we partition the set of marked words into three classes:
\begin{enumerate}
\item The empty word, which has weight $q^0$.
\item Those words ending with a letter which is not part of a marked descent.
\item Those words ending with a marked descent.
\end{enumerate}
If a marked word ends with a marked descent, it must end in a marked descending run $w_i > w_{i+1} > \cdots > w_\ell$, where each component descent is marked.  We remove the maximal terminal marked run (also called a \emph{cluster}), leaving a word in the second category above.  To continue the example above, $(2637541, \{2,5,6\})$ has $541$ as its maximal terminal marked run.  To contrast, note that $(2637541, \{2,4,6\})$ has only $41$ as its maximal terminal marked run, since $54$ is not marked.  Since weights are multiplicative, the weight of $(2637541, \{2,5,6\})$ is the product of the weights of $(2637541, \{5,6\})$ (the cluster) and $(2637, \{2\})$ (the rest).  Furthermore, rearranging the letters outside the cluster does not affect the weight of the cluster, e.g. $(2637541, \{5,6\})$ and $(7362541, \{5,6\})$ have equal weight.  Thus we may fix a cluster and let the prefix vary, guiding our next transformations.

For any subset $T\subseteq [n]$, identify $T={t_1>t_2>\cdots>t_k}$ with its corresponding descending run $t_1 t_2 \cdots t_{k}$, and so write $wT=wt_1 t_2 \cdots t_k$.  Algebraically, the above partition translates into:
\begin{eqnarray*}
  F_{\sigma}(q,x)&=& \sum_{w\in A^*} x^w \sum_{S\subseteq Des(w)} \prod_{i\in S} (q^{\sigma_d(w,i)}-1) \\
        &=& 1 + \sum_{i\in [n]} \sum_{w\in A^*} x^{wi} \sum_{S\subseteq Des(w)} \prod_{i\in S} (q^{\sigma_d(w,i)}-1)\\
        && \hspace{.2 in} +\sum_{T\subseteq [n], |T|\geq 2} \sum_{w\in A^*} x^{wT} \sum_{S\subseteq Des(w)} \prod_{i\in S} (q^{\sigma_d(w,i)}-1) \left( \prod_{j=1}^{|T|-1} (q^{\sigma_d(wT,\ell(w)+j)}-1)\right) \\
        &=&1+ \sum_{T\subseteq [n],T\neq \emptyset} \sum_{w\in A^*} x^{wT} \left(\prod_{j=1}^{|T|-1} (q^{\sigma_d(wT,\ell(w)+j)}-1)\right) \sum_{S\subseteq Des(w)} \prod_{i\in S} (q^{\sigma_d(w,i)}-1) \\
\end{eqnarray*}

To rephrase the preceding observation, $\sigma_d(wT,\ell(w)+j)$ is dependent only on the values of the $\nu_i(w)$, and not on the order of the letters in $w$ itself.  For example, it is clear to see that $(a-cb)_d(wT,\ell(w)+j)=\nu_1(w)+\nu_2(w)+\cdots+\nu_{t_{j+1}-1}(w)$.  Therefore, let $m=(m_1,\ldots,m_n)\in \mathbb{N}^n$ and partition $A^*$ into classes $A^*_m = \{ w\in A^* : \nu_i(w)=m_i \text{ for all } i\}$.  Note $A^*_m$ is the set of multiset permutations with $m$ giving the multiplicities of each element.  Let $\tilde{\sigma}_d(m,T,j)=\sigma_d(wT,m_1+\cdots m_n + j)$, where $w\in A^*_m$ (it does not matter which $w$ you choose).  This is the descent-weight of the $j^{th}$ descent in $T$ in the word $wT$.  For example, $\widetilde{(a-cb)}_d(m,T,j)=m_1+m_2+\cdots+m_{t_{j+1}-1}$.
 Hence we can write $F_\sigma(q,x)$ as
\[
F_{\sigma}(q,x) = 1+ \sum_{m\in \mathbb{N}^n} x_1^{m_1} \cdots x_n^{m_n} \sum_{T\subseteq [n], T \neq \emptyset} x^T \prod_{j=1}^{|T|-1}(q^{\tilde{\sigma}_d(m,T,j)}-1) \sum_{w\in A^*_m} \sum_{S\subseteq Des(w)} \prod_{i\in S} (q^{\sigma_d(w,i)}-1)
\]

Let us focus on $\sum_{w\in A^*_m} q^{\sigma(w)}$ (the coefficient of $x_1^{m_1}\cdots x_n^{m_n}$ in $F_\sigma(q,x)$), which we denote $F_\sigma(m,q)$.  We will write recurrences for $F_\sigma(m,q)$ in the form of operators $P_\sigma$ so that $P_\sigma F_\sigma = F_\sigma$.  Let $E_i^{-1}$ be the shift operator for $m_i$, that is, \[E_i^{-1} f(m_1, \ldots, m_i, \ldots, m_n) = f(m_1, \ldots, m_i-1, \ldots, m_n)\] and $E_T^{-1} = \prod_{i\in T} E_i^{-1}$.  Then it follows that 
\[
F_\sigma(m,q) = \sum_{T\subseteq [n], T\neq \emptyset} E_T^{-1} \prod_{j=1}^{|T|-1} (q^{\tilde{\sigma}_d(m,T,j)}-1) F_\sigma(m,q)
\]
for $m=(m_1, \ldots, m_n)$ and each $m_i>0$.  In the case that some $m_i=0$, then $F_\sigma(m,q) = F_\sigma(\hat{m},q)$ for $\hat{m}=(m_1, \ldots, m_{i-1},m_{i+1}, \ldots, m_n)$, since we may reduce by one each letter greater than $i$.  Furthermore, if $m_i=0$ for all $i$, then $F_\sigma(m,q)=1$.  These provide the initial conditions for the recurrence above.

As we apply each $E_T^{-1}$ operator, we will need to consider the shift's effect on $\tilde{\sigma}_d(m,T,j)$.  It is easily checked that 
\begin{eqnarray*}
E_T^{-1} \widetilde{(a-cb)}_d (m,T,j) &=& m_1+m_2+\cdots+m_{t_{j+1}-1}-|T|+j+1 \\
   &=&\widetilde{(a-cb)}_d (m,T,j) - |T|+j+1
\end{eqnarray*}
since the descent $t_j t_{j+1}$ follows all occurences of letters less than $t_{j+1}$ except the $|T|-(j+1)$ which lie in $T$ itself.
  Similarly, 
\[E_T^{-1} \widetilde{(a-ba)}_d(m,T,j) = \widetilde{(a-ba)}_d(m,T,j) - 1\]
\[E_T^{-1} \widetilde{(b-ba)}_d(m,T,j) = \widetilde{(b-ba)}_d(m,T,j) - 1.\]
For $\sigma \in \{(ba), (b-ca), (c-ba)\}$, $E_T^{-1} \tilde{\sigma}_d(m,T,j) = \tilde{\sigma}_d(m,T,j)$.  The operators for each descent-based $\sigma$ are summarized in the table below.

\begin{table}[ht]
	\title{Operators $P_\sigma$ for descent-based $\sigma$}
		{\begin{tabular}{|c|c|}
		  \hline
			$\sigma$ & $P_\sigma$ \\
			\hline
			$(ba)  $ & $\sum_{T\subseteq [n], T\neq \emptyset}  \prod_{j=1}^{|T|-1} (q^{1}-1) E_T^{-1} $ \\ 
			$(a-cb)$ & $\sum_{T\subseteq [n], T\neq \emptyset}  \prod_{j=1}^{|T|-1} (q^{m_1+m_2+\cdots+m_{t_{j+1}-1}-|T|+j+1}-1) E_T^{-1} $ \\ 
			$(b-ca)$ & $\sum_{T\subseteq [n], T\neq \emptyset}  \prod_{j=1}^{|T|-1} (q^{m_{t_{j+1}+1}+m_{t_{j+1}+2}+\cdots+m_{t_{j}-1}}-1) E_T^{-1} $ \\ 
			$(c-ba)$ & $\sum_{T\subseteq [n], T\neq \emptyset}  \prod_{j=1}^{|T|-1} (q^{m_{t_{j}+1}+m_{t_{j}+2}+\cdots+m_{n}}-1) E_T^{-1} $ \\ 
			$(a-ba)$ & $\sum_{T\subseteq [n], T\neq \emptyset}  \prod_{j=1}^{|T|-1} (q^{m_{t_{j+1}}-1}-1) E_T^{-1} $ \\ 
			$(b-ba)$ & $\sum_{T\subseteq [n], T\neq \emptyset}  \prod_{j=1}^{|T|-1} (q^{m_{t_{j}}-1}-1) E_T^{-1} $ \\ 
			\hline
		\end{tabular}}
	\label{Des-based Operators}
\end{table}

For rise-based statistics, the process remains unchanged except now one considers marked rises (i.e. subsets of $Rise(w)=\{i: w_i<w_{i+1}\}$).  To keep the relative ordering fo the $t_i\in T$ consistent, we consider what happens upon the removal of the ascending run $T^r$ from $w\!T^r$.  
\[F_\sigma(q,x)=1+ \sum_{T\subseteq [n],T\neq \emptyset} \sum_{w\in A^*} x^{wT^r} \left(\prod_{j=1}^{|T|-1} (q^{\sigma_d(wT^r,\ell(w)+j)}-1)\right).\]  Following the same manipulations as for descent=based $\sigma$, we obtain the corresponding operators shown in the table below.

\begin{table}[ht]
	\title{Operators $P_\sigma$ for rise-based $\sigma$}
		{\begin{tabular}{|c|c|}
		  \hline
			$\sigma$ & $P_\sigma$ \\
			\hline
			$(ab)  $ & $\sum_{T\subseteq [n], T\neq \emptyset}  \prod_{j=1}^{|T|-1} (q^{1}-1) E_T^{-1} $ \\ 
			$(a-bc)$ & $\sum_{T\subseteq [n], T\neq \emptyset}  \prod_{j=1}^{|T|-1} (q^{m_1+m_2+\cdots+m_{t_{j+1}-1}}-1) E_T^{-1} $ \\ 
			$(b-ac)$ & $\sum_{T\subseteq [n], T\neq \emptyset}  \prod_{j=1}^{|T|-1} (q^{m_{t_{j+1}+1}+m_{t_{j+1}+2}+\cdots+m_{t_{j}-1}}-1) E_T^{-1} $ \\
			$(c-ab)$ & $\sum_{T\subseteq [n], T\neq \emptyset}  \prod_{j=1}^{|T|-1} (q^{m_{t_{j}+1}+m_{t_{j}+2}+\cdots+m_{n}-j+1}-1) E_T^{-1} $ \\ 
			$(a-ab)$ & $\sum_{T\subseteq [n], T\neq \emptyset}  \prod_{j=1}^{|T|-1} (q^{m_{t_{j+1}}-1}-1) E_T^{-1} $ \\ 
			$(b-ab)$ & $\sum_{T\subseteq [n], T\neq \emptyset}  \prod_{j=1}^{|T|-1} (q^{m_{t_{j}}-1}-1) E_T^{-1} $ \\ 		
			\hline
	 \end{tabular}}
	\label{Rise-based Operators}
\end{table}

Since the reversal operator is a bijection in $A^*_m$, we see $P_{\sigma^r}=P_{\sigma}$.  Hence we have recurrences for $F_\sigma(m,q)$ for each $\sigma$ of type-$(2,1)$ as well.  Note that here the clusters are initial marked runs instead of terminal marked runs.

Some operators appear multiple times among Tables \ref{Des-based Operators} and \ref{Rise-based Operators}.  Since all patterns have the same initial conditions, we have proven the following equidistribution results.

\begin{theorem}
The patterns $(b-ac)$ and $(b-ca)$ are equidistributed over $A^*_m$ for any alphabet vector $\mathbf{m}$.  Similarly $(a-ab)$ and $(a-ba)$ are equidistributed over $A^*_m$, and so are $(b-ab)$ and $(b-ba)$.
\end{theorem}

These are the only such non-trivial equi-distribution classes over multiset permutations among the patterns listed above, as evidenced by computing the distribution over $A^*_{(1,1,1,2)}$ (see table below).

\begin{table}[ht]
	\title{Distributions of descent-based statistics over $A^*_{(1,1,1,2)}$}
		{\begin{tabular}{|c|c|}
		 \hline
			$\sigma $   & $F_\sigma(q,(1,1,1,2))$ \\
		 \hline
			$(b a)    $ & $1+18q+33q^2+8q^3$ \\
			$(a - c b)$ & $31+17q+11q^2+q^3$ \\
			$(b - c a)$ & $28+23q+8q^2+q^3$ \\
			$(c - b a)$ & $37+10q+9q^2+3q^3+q^5$ \\
			$(a - b a)$ & $60$ \\
			$(b - b a)$ & $24+36q$ \\
			$(a b)    $ & $1+18q+33q^2+8q^3$ \\
			$(a - b c)$ & $31+20q+5q^2+4q^3$ \\
			$(b - a c)$ & $28+23q+8q^2+q^3$ \\
			$(c - a b)$ & $37+9q+10q^2+3q^3+q^4$ \\
			$(a - a b)$ & $60$ \\
			$(b - a b)$ & $24+36q$ \\
		 \hline
		\end{tabular}}
	\label{Distributions}
\end{table}

This method's strength lies in its applicability to multistatistics.  Let $\sigma$ and $\tau$ be two descent-based pattern functions and let $F_{(\sigma,\tau)} (m,q,t):=\sum_{w\in A^*_m} q^{\sigma(w)} t^{\tau(w)}$.  Then just as before we may write $F_{\sigma,\tau}$ as 
\[
  F_{\sigma,\tau}(m,q,t)= \sum_{w\in A^*_m} \sum_{S\subseteq Des(w)} \prod_{i\in S} (q^{\sigma_d(w,i)} t^{\tau_d(w,i)}-1)   
\]

All of the above transformations hold, leading to the operator
\[P_{\sigma,\tau} = \sum_{T\subseteq [n], T\neq \emptyset} E_T^{-1} \prod_{j=1}^{|T|-1} (q^{\tilde{\sigma}_d(m,T,j)} t^{\tilde{\tau}_d(m,T,j)} -1).\]

Of course this generalizes in the same way to multistatistics with any number of descent-based patterns of type-(1,2).  Furthermore, specializing indeterminates (except $x$) allows one to consider compound statistics, such as the major index $maj=(ba)+(a-cb)+(b-ca)+(c-ba)+(b-ba)+(a-ba)$ or multistatistic distributions like the Euler-Mahonian polynomials $\sum_{\pi \in S_n} t^{maj(\pi)} q^{(ba)(\pi)} $.  Likewise one may consider multistatistics involving rise-based pattern functions.  The greatest restriction is that one may only combine descent-based statistics, or combine rise-based statistics, and must keep within type-(1,2) or type-(2,1).  In short, one may combine patterns from only one of the four boxes below.

\begin{center}
	\begin{tabular}{|c|c|c|}
		 \hline
		   $\begin{array}{c} (ba), (a-ba), (b-ba)\\ (a-cb), (b-ca), (c-ba) \end{array}$ & $\begin{array}{c} (ab), (a-ab), (b-ab) \\ (a-bc), (b-ac), (c-ab) \end{array}$ \\
		 \hline
		   $\begin{array}{c} (ab), (ab-a), (ab-b)\\ (bc-a), (ac-b), (ab-c) \end{array}$ & $\begin{array}{c} (ba), (ba-a), (ba-a) \\ (cb-a), (ca-b), (ba-c) \end{array}$ \\
		 \hline
	\end{tabular}
\end{center}

Problems combining descent-based and rise-based patterns arise since both descents and rises are marked, resulting in non-monotone clusters.  Thus for each $T\subseteq [n]$ there is more than one cluster involving the letters in $T$.  Furthermore one must consider submultisets $T$ to compute distributions over multiset permutations.  A variant of the above methods may still work for such multistatistics if one sacrifices speed.  As mentioned previously, the clusters for type-(1,2) patterns appear at the end of the word whereas clusters for type-(2,1) appear at the beginning.  How to rectify this disparity eludes the author.  

\section{Distributions over $S_n$}

We now move to the special case of permutations, i.e. computing $F_{\sigma}(m)$ for $m=(1,1,\ldots, 1)\in \mathbb{N}^n$, which we will denote $F_{\sigma}(n)$.  In this case it is clear that the shift operator $E^{-1}_T$ translates to the shift $N^{-|T|} F(n)=F(n-|T|)$.  For each of $(ba)$, $(a-cb)$, $(b-ca)$, and $(c-ba)$ we will determine a recurrence to recursively compute \[F_{\sigma}(n)=\sum_{\pi \in S_n} q^{\sigma(\pi)}.\]  These recurrences allow for polynomial-time computation, and the same method should allow for polynomial-time computation for multistatistics following the restrictions discussed at the end of the previous section.  In each subsection below, the $\sigma$ subscript for $F$ is omitted.  Of course for each pattern $\sigma$ discussed below, the same recurrences will hold for the distributions of $\sigma^r$, $\sigma^c$, and $\sigma^{c\;r}$.

\subsection{The pattern (ba)}
   From table \ref{Des-based Operators} we get the recurrence
    \begin{eqnarray*}
     F(n) &=& \sum_{T\subseteq [n], T\neq \emptyset}  \prod_{j=1}^{|T|-1} (q^{1}-1) F(n-|T|) \\
          &=& \sum_{k=1}^n \binom{n}{k} (q-1)^{k-1} F(n-k)
    \end{eqnarray*}
   Observe that $q F(n)$ are the well-known Eulerian polynomials.

\subsection{The pattern (a-cb)}
    From table \ref{Des-based Operators} we get the recurrence
    \[F(n) = \sum_{k=1}^n \sum_{n \geq t_1 > \cdots > t_k \geq 1}  \prod_{j=1}^{k-1} (q^{t_{j+1}+j-k}-1) F(n-k)\]
  since $m_1+m_2+\cdots+m_{t_{j+1}-1} = t_{j+1}-1$.  Define 
   \[a(n,k)=\sum_{n \geq t_1 > \cdots > t_k \geq 1}  \prod_{j=1}^{k-1} (q^{t_{j+1}+j-k}-1)\]
  so that $F(n)=\sum a(n,k) F(n-k)$, and let 
    \[b(n,k)=\sum_{n \geq t_1 > \cdots > t_k \geq 1}  \prod_{j=1}^{k} (q^{t_{j}+j-k-1}-1).\]
    Conditioning on whether $t_1=n$, we get that 
    \begin{eqnarray*}
      a(n,k)&=& \sum_{n-1 \geq t_1 > \cdots > t_k \geq 1}  \prod_{j=1}^{k-1} (q^{t_{j+1}+j-k}-1) + \sum_{n-1 \geq t_2 > \cdots > t_k \geq 1}  \prod_{j=1}^{k-1} (q^{t_{j+1}+j-k}-1) \\
        &=& a(n-1,k) + \sum_{n-1 \geq t_1 > \cdots > t_{k-1} \geq 1}  \prod_{j=1}^{k-1} (q^{t_{j}+j-k}-1)\\
        &=& a(n-1,k) + b(n-1,k-1)
     \end{eqnarray*}
     Similarly we can derive a recurrence for $b(n,k)$ by conditioning on $t_1$.
     \begin{eqnarray*}
      b(n,k)&=&\sum_{n \geq t_1 > \cdots > t_k \geq 1}  \prod_{j=1}^{k} (q^{t_{j}+j-k-1}-1) \\
          &=& \sum_{n-1 \geq t_1 > \cdots > t_k \geq 1}  \prod_{j=1}^{k} (q^{t_{j}+j-k-1}-1) +\sum_{n-1 \geq t_2 > \cdots > t_k \geq 1}  (q^{n-k}-1)  \prod_{j=2}^{k} (q^{t_{j}+j-k-1}-1) \\
           &=&b(n-1,k)+(q^{n-k}-1) \sum_{n-1 \geq t_1 > \cdots > t_{k-1} \geq 1}  \prod_{j=1}^{k-1} (q^{t_{j}+j-k}-1) \\
           &=&b(n-1,k)+(q^{n-k}-1) b(n-1,k-1)
      \end{eqnarray*}
     These recurrences, along with the initial conditions below yield fast computation of $F(n)$.
\begin{center}
	\begin{tabular}{ll}
		$a(n,1)=n$               & $b(n,1) = \sum_{i=1}^n (q^{i-1}-1)$ \\
		$a(n,k)=0$ for $k\geq n$ & $b(n,k)=0$ for $k\geq n$ 
	\end{tabular}
\end{center}
  This recurrence confirms the values listed in Table 2 in \cite{ClaesMans}.

\subsection{The pattern (b-ca)}
  As in the case for $(a-cb)$, we get a recurrence of the form $F(n)=\sum_{k=1}^n a(n,k) F(n-k)$, where now we get
  \[a(n,k)=\sum_{n \geq t_1 > \cdots > t_k \geq 1}  \prod_{j=1}^{k-1} (q^{t_{j}-t_{j+1}-1}-1).\]
  We will also make use of a secondary function $b(n,k)$, defined as
  \[b(n,k)=\sum_{n \geq t_1 > \cdots > t_k \geq 1}  q^{-t_1} \prod_{j=1}^{k-1} (q^{t_{j}-t_{j+1}-1}-1).\]
  Again conditioning on whether $t_1=n$, we obtain
  \begin{eqnarray*}
    a(n,k) &=& a(n-1,k) + \left( \sum_{n-1 \geq t_2 > \cdots > t_k \geq 1}  (q^{n-t_2-1}-1)\prod_{j=1}^{k-1} (q^{t_{j}-t_{j+1}-1}-1) \right) \\
     &=& a(n-1,k) + q^{n-1} \sum_{n-1 \geq t_2 > \cdots > t_k \geq 1} q^{-t_2} \prod_{j=2}^{k-1} (q^{t_{j}-t_{j+1}-1}-1) - \sum_{n-1 \geq t_2 > \cdots > t_k \geq 1}  \prod_{j=1}^{k-1} (q^{t_{j}-t_{j+1}-1}-1) \\
     &=& a(n-1,k) + q^{n-1} b(n-1,k-1) - a(n-1,k-1) \\
   && \\
    b(n,k) &=& b(n-1,k) + \sum_{n-1 \geq t_2 > \cdots > t_k \geq 1}  q^{-n} (q^{n-t_2-1}-1) \prod_{j=2}^{k-1} (q^{t_{j}-t_{j+1}-1}-1) \\
      &=&b(n-1,k) + q^{-n} \left( \sum_{n-1 \geq t_2 > \cdots > t_k \geq 1}  (q^{n-t_2-1}-1)\prod_{j=1}^{k-1} (q^{t_{j}-t_{j+1}-1}-1) \right) \\
      &=&b(n-1,k) + q^{-n} \left( q^{n-1} b(n-1,k-1) - a(n-1,k-1) \right)  
  \end{eqnarray*}
  With the initial conditions below we can quickly compute terms of $F(n)$.
 \begin{center}
	\begin{tabular}{ll}
		$a(n,1)=n$               & $b(n,1) = \sum_{i=1}^n (q^{-i}-1)$ \\
		$a(n,k)=0$ for $k\geq n$ & $b(n,k)=0$ for $k\geq n$ 
	\end{tabular}
\end{center}

Note that Parviainen\cite{Parv} has found closed-form formulas for the coefficient of $q^k$ in $F(n)$ for $1\leq k\leq8$ and provides the automated method for higher $k$.  This recurrence confirms the values listed in Table 3 in \cite{ClaesMans}.

\subsection{The pattern (c-ba)}
  We again seek a recurrence of the form $F(n)=\sum_{k=1}^n a(n,k) F(n-k)$, where 
  \[a(n,k) = \sum_{n \geq t_1 > \cdots > t_k \geq 1}  \prod_{j=1}^{k-1} (q^{n-t_{j}}-1).\]
  We will also need to use the secondary functions
  \[b(n,k) = \sum_{n \geq t_1 > \cdots > t_k \geq 1} t_k \prod_{j=1}^{k} (q^{n-t_{j}}-1),\]
  and 
  \[c(n,k) = \sum_{n \geq t_1 > \cdots > t_k \geq 1}  \prod_{j=1}^{k} (q^{n-t_{j}}-1).\]
  First observe that $a(n,k)=b(n-1,k-1)$, since
  \begin{eqnarray*}
    a(n,k) &=& \sum_{n \geq t_1 > \cdots > t_{k-1} \geq 2} (t_k-1) \prod_{j=1}^{k-1} (q^{n-t_{j}}-1) \\
      &=& \sum_{n-1 \geq t_1 > \cdots > t_{k-1} \geq 1} t_k \prod_{j=1}^{k-1} (q^{n-1-t_{j}}-1) \\
      &=& b(n-1,k-1)
  \end{eqnarray*}
  Now consider the sum in $b(n,k)$ and condition on whether $t_k=1$.  Then we get
  \begin{eqnarray*}
    b(n,k) &=& \sum_{n \geq t_1 > \cdots > t_k \geq 2} t_k \prod_{j=1}^{k} (q^{n-t_{j}}-1) + \sum_{n \geq t_1 > \cdots > t_{k-1} \geq 2, t_k=1}  \prod_{j=1}^{k} (q^{n-t_{j}}-1)  \\
     &=& \sum_{n-1 \geq t_1 > \cdots > t_k \geq 1} (t_k+1) \prod_{j=1}^{k} (q^{n-1-t_{j}}-1) + (q^{n-1}-1)\sum_{n-1 \geq t_1 > \cdots > t_{k-1} \geq 1, t_k=1}  \prod_{j=1}^{k-1} (q^{n-1-t_{j}}-1)  \\
     &=& b(n-1,k)+c(n-1,k) + (q^{n-1}-1) c(n-1,k-1)
   \end{eqnarray*}
  We also condition on $t_k$ to get a recurrence for $c(n,k)$.
  \begin{eqnarray*}
    c(n,k) &=& \sum_{n \geq t_1 > \cdots > t_k \geq 2}  \prod_{j=1}^{k} (q^{n-t_{j}}-1) + \sum_{n \geq t_1 > \cdots > t_{k-1} \geq 2, t_k=1} (q^{n-1}-1) \prod_{j=1}^{k-1} (q^{n-t_{j}}-1) \\
           &=& \sum_{n-1 \geq t_1 > \cdots > t_k \geq 1}  \prod_{j=1}^{k} (q^{n-1-t_{j}}-1) + (q^{n-1}-1) \sum_{n-1 \geq t_1 > \cdots > t_{k-1} \geq 1}  \prod_{j=1}^{k-1} (q^{n-1-t_{j}}-1) \\
           &=& c(n-1,k) + (q^{n-1}-1) c(n-1,k-1)
   \end{eqnarray*}

This recurrence confirms the values listed in Table 1 in \cite{ClaesMans}.

\section{Accompanying Maple Package \textit{ClusterGPP}}

A corresponding package of Maple procedures illustrating the above methods can be downloaded at this paper's website:
\[  \underline{ http://math.rutgers.edu/ \tilde{\;}baxter/ClusterGPP/index.html} \]
Load the file into Maple using the \textit{read} command.  The built-in help command, \textit{Help()} will outline the program for you, describing the syntax for each of the main procedures (listed for your convenience in the table below.

The reader is encouraged to experiment with the above package to see resulting distributions.

\begin{table}[ht]
	\title{Procedures in the Maple package \textit{ClusterGPP}}
		{\begin{tabular}{|l|p{3.5 in}|}
		 \hline
			Procedure Name        & Description \\
		\hline
			\textit{PatternCount} & Counts the number of occurences of a given pattern in a given word \\
		\hline
			\textit{BFdist} & Computes the distribution of a given set of patterns over a given set of words. \\
		\hline	
			\textit{Rdist} & Computes the distribution $F_{\sigma}(m,q)$ for a single pattern $\sigma$ using the operator $P_\sigma$. \\
		\hline
			\textit{RdistM} & Recursively computes the distribution of multiple patterns, e.g. $F_{\sigma,\tau}(m)$ using the operator $P_{\sigma,\tau}$.  This subsumes \textit{Rdist}, as singletons may be entered. \\
   \hline
		\end{tabular}}
		\label{Procedures}
\end{table}

\section{Conclusions and Future Directions}

The methods outlined above yield recurrences for the distribution of any number of multistatistics generated by certain combinations of generalized permutation patterns.  Besides being interesting in their own right, setting the indeterminates to 0 in the distributions gives the corresponding pattern avoidance results, such as those results in \cite{Claes01}.  For example, setting $q=0$ in the recurrence for the pattern (a-cb) above gives an alternate recurrence for the generation of the Bell numbers.  Questions of pattern packing, as studied in \cite{PattPack1, PattPack2, PattPack3}, ask for the highest number of occurences of a given pattern or patterns.  This is simply the degree of the distributions we have calculated.  For example, the first 20 terms of the sequence \[a(n):=\max \{(a-cb)(\pi) : \pi \in S_n \}\] can be computed in under a minute using procedures from \textit{ClusterGPP}.  The interested reader may adapt the above procedures to generate recurrences for the degrees of the polynomials.

The methods above would need massaging to extend them to multistatistics involving both rise- and descent-based patterns, for example the multi-pattern function $(a-cb)+(a-bc)$.  One must overcome the hurdle of weight-counting both descents and rises simultaneously, and determining what to do with runs of marked descents and rises.  To extend the methods to 4-patterns, such as $(a-c-db)$, the weights on the descents would get much more complicated than they are above.  Of course, certain combinations of 4-patterns are already amenable to our methods, such as $(a-c-db)+(c-a-db) = (b-ca)+(a-cb)$.

\end{document}